\documentclass[12pt]{amsart}

\usepackage{amsthm}

\usepackage{amssymb}
\usepackage{verbatim}
\usepackage[curve,matrix,arrow]{xy}
\usepackage{amsmath,amsfonts}
\usepackage{graphicx}
\usepackage{epsf}

\newcommand{\thmref}[1]{Theorem~\ref{#1}}
\newcommand{\propref}[1]{Proposition~\ref{#1}}
\newcommand{\lemref}[1]{Lemma~\ref{#1}}
\newcommand{\eqnref}[1]{equation~(\ref{#1})}

\newcommand{\corref}[1]{Corollary~\ref{#1}}

  {\end{list}}

\def\li{L_{i}}

\def\RR{{\mathbb R}}

\newtheorem{theorem}{Theorem}[section]

\newtheorem{corollary}[theorem]{Corollary}

\newtheorem{lemma}[theorem]{Lemma}
\newtheorem{proposition}[theorem]{Proposition}

\newtheorem{remark}[theorem]{Remark}

\newcommand{\jj}[3]{j_{#1}(#2,#3)}

\newcommand{\ga}{\Gamma}

\newcommand{\ee}[1]{E(#1)}
\newcommand{\vv}[1]{V(#1)}
\newcommand{\va}{\upsilon}

\newcommand{\p}[1]{ \left( #1 \right) }

\def\<{\langle }
\def\>{\rangle }

\def\tr{\mathrm{tr}}
\def\diag{\text{diag}}
\newcommand{\am}{\mathrm{A}}
\newcommand{\ay}{\mathrm{Y}}

\newcommand{\dm}{\mathrm{D}}
\newcommand{\mm}{\mathrm{M}}

\newcommand{\idm}{\mathrm{I}}
\newcommand{\jm}{\mathrm{J}}
\newcommand{\hm}{\mathrm{H}}
\newcommand{\om}{\mathrm{O}}

\newcommand{\lm}{\mathrm{L}}
\newcommand{\plm}{\mathrm{L^+}}

\newcommand{\lpq}{l_{pq}}

\begin{document}

\title[Generalized Foster's Identities]
{Generalized Foster's Identities}


\author{Zubeyir Cinkir}
\address{Zubeyir Cinkir\\
Department of Mathematics\\
University of Georgia\\
Athens, Georgia 30602\\
USA}
\email{cinkir@math.uga.edu}

\thanks{I would like to thank Dr. Robert Rumely for his continued support and the discussions about this paper}

\begin{abstract}
Foster's network theorems and their extensions to higher orders involve resistance values and conductances. We establish identities concerning voltage values and conductances. Our identities are analogous to the extended Foster's identities.
\end{abstract}

\maketitle

\section{Introduction}\label{sec introduction}

Foster's first identity, which is also known as Foster's network theorem, was initially proved by R. M. Foster \cite{Fo1} in $1949$ by using Kirchhoff's rule.
In $1990$, D. Copersmith, P. Doyle, P. Raghavan, and M. Snir \cite{CDRS} gave another proof of
this theorem by using random walks on graphs.
Independently, in $1991$, P. Tetali \cite{TP1} gave a third proof of this
theorem again by using random walks on graphs. (See
also \cite[Theorem 25 and exercise 23 in Chapter IX]{BB}). Later in
$1994$, Tetali \cite{TP2} gave a fourth probabilistic proof by
using an elementary identity for ergodic Markov chains.

In the late $1980$'s, T. Chinburg and R. Rumely
introduced a canonical measure $\mu_{can}$ on a metrized graph.
As a consequence of the fact that $\mu_{can}$ has total mass $1$,
Chinburg and Rumely \cite[remark on pg. 26]{CR} obtained an identity.
In $2003$, M. Baker and X. Faber \cite[Corollary 6]{BF} showed
that this identity is equivalent to the Foster's first identity.
This became another interesting proof of the Foster's first identity.

In another direction, L. Weinberg \cite{W} in $1958$, D. J. Klein and M. Randi\'c \cite[Corollary C1]{D-M} in $1993$, and R. Bapat \cite[Lemma $2$]{RB1} in $2004$ proved identities that are equivalent to the Foster's first identity. The properties of discrete Laplacian and its pseudo inverse were used in the articles \cite{D-M} and \cite{RB1}.

In $1961$, Foster \cite{Fo2} proved another identity which we call Foster's
second identity. In 2002, J. L. Palacios \cite{P} gave
a probabilistic proof of the Foster's second identity,
generalizing the arguments used in \cite{TP2}. He also showed how to
extend the Foster's identities to graph paths with more edges and he proved a third identity.
The Foster's first and second identities, and the third identity due to
Palacios are about the sums of the resistance values over the graph
paths consisting of one, two and three edges, respectively.

In $2007$, we \cite[Section 5.5]{C1} proved four identities that involve voltage values rather than resistance values.
Our method was to use the properties of discrete laplacian and its pseudo inverse. As a corollary to these identities,
we obtained Foster's 1st, 2nd, 3rd, and 4th identities.

In $2008$, E. Bendito, A. Carmona, A. M. Encinas and J. M. Gesto
\cite[Proposition 2.3]{BCEG} extended the Foster's identities to higher orders. Their formula contains all previously known Foster's identities as specific cases.

In this paper, we extend the identities concerning voltage values \cite[Section 5.5]{C1} to higher orders. This can be found in \thmref{thm main}, which is our main result.
In this way, we generalize the extended Foster's identities.
We provide two proofs for \thmref{thm main}. Both of the proofs rely on the properties of the discrete Laplacian and its pseudo inverse.
In the first proof, we use a remarkable relation between the resistance values and the pseudo inverse $\plm$ of the discrete Laplacian \cite{RB1} \cite{RB2} (see also \lemref{lem disc2}) to express the voltage values by using $\plm$ (see \lemref{lem voltage1}).
Similarly, in the second proof, we use an equally interesting relation between the resistance values and some ``equilibrium measures'' $\nu^i$'s \cite[Proposition 2.1]{BCEG} to express the voltage values by using $\nu^i$'s (see \propref{propeqn voltage}).
More information about the equilibrium measure can be found in the article \cite{BCEG} and the related references therein.


Note that there is a  $1-1$ correspondence between the metrized graphs and the
equivalence class of finite connected weighted graphs in which the weight of an edge is
the reciprocal of its length \cite[Lemma 2.2]{BRh}. We can also view such a graph as a
resistive electrical network in which the resistance of each edge is
the same as its length. We will work with metrized graphs in this paper. However, the results and their proofs are
valid for the corresponding weighted graphs and the resistive electrical networks.

\section{The voltage and resistance functions}\label{sec resistances}

A \textit{metrized graph} $\ga$ is a finite connected graph equipped with a distinguished
parametrization of each of its edges. One can find other definitions of metrized graphs in the articles \cite{CR}, \cite{BRh},
\cite[Appendix]{Zh1} and \cite{BF}.

A metrized graph can have multiple edges and self-loops. For any given $p \in \ga$,
the number of directions
emanating from $p$ will be called the \textit{valence} of $p$, and will be denoted by
$\va(p)$. By definition, there can be only finitely many $p \in \ga$ with $\va(p)\not=2$.

For a metrized graph $\ga$, we will denote its set of vertices by $\vv{\ga}$.
We require that $\vv{\ga}$ be finite and non-empty and that $p \in \vv{\ga}$ for each $p \in \ga$ if
$\va(p)\not=2$. For a given metrized graph $\ga$, it is possible to enlarge the
vertex set $\vv{\ga}$ by considering more additional points of valence $2$ as vertices.

For a given metrized graph $\ga$, the set of edges of $\ga$ is the set of closed line segments with end points
in  $\vv{\ga}$. We will denote the set of edges of $\ga$ by $\ee{\ga}$.
We will denote $\# (\vv{\ga})$ and $\# (\ee{\ga})$ by $n$ and $e$, respectively if there is no danger of confusion.
We denote the length of an edge $e_i \in \ee{\ga}$ by $\li$.

In the article \cite{CR}, the \textit{voltage function} $j_{z}(x,y)$
is defined and studied as a function of $x, y, z \in \Gamma$. For
fixed $z$ and $y$ it has the following physical interpretation: when
$\Gamma$ is viewed as a resistive electric circuit with terminals at
$z$ and $y$, with the resistance in each edge given by its length,
then $j_{z}(x,y)$ is the voltage difference between $x$ and $z$,
when unit current enters at $y$ and exits at $z$ (with reference
voltage 0 at $z$).

For any $x$, $y$, $z$ in $\ga$, the voltage function $j_x(y,z)$ on
$\ga$ is a symmetric function in $y$ and $z$, i.e., $j_x(y,z)=j_x(z,y)$. Moreover, for any $x$, $y$, $z$ in $\ga$ it satisfies
$j_x(x,z)=0$ and $j_x(y,y)=r(x,y)$, where $r(x,y)$ is the \textit{resistance
function} on $\ga$ (for more information see the articles \cite{CR} and \cite[sec 1.5 and sec 6]{BRh}, and \cite[Appendix]{Zh1}).

\begin{figure}
\centering
\includegraphics[scale=0.7]{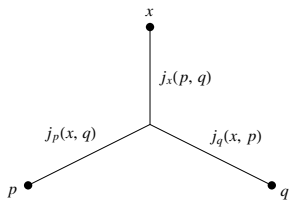} \caption{Circuit reduction with reference to $3$ points $x$, $p$ and $q$.} \label{fig xpq1new2}
\end{figure}

For any $x \in \ga$, by circuit theory we can transform $\ga$ to an
$Y$-shaped graph with the same resistances between $x$, $p$, and
$q$ as in $\ga$ (see \cite[Section 2]{C2} for more details).
This is shown in Figure \ref{fig xpq1new2}, with the corresponding
voltage values on each segment.
Therefore,
\begin{equation}\label{eqn1.1}
\begin{split}
\\r(p,x) = \jj{p}{x}{q}+\jj{x}{p}{q}, \; \, r(q,x) = \jj{q}{x}{p}+\jj{x}{p}{q},
\; \, r(p,q) = \jj{q}{x}{p}+\jj{p}{x}{q}.
\end{split}
\end{equation}
Then it follows from (\ref{eqn1.1}) that
\begin{equation}\label{eqn1.1b}
\begin{split}
2j_p(x,q)=r(p,x)+r(p,q)-r(q,x), \qquad \text{for any $p$, $q$, $x$ in $\ga$}.
\end{split}
\end{equation}

\section{The discrete Laplacian $\lm$ and the pseudo inverse $\plm$}
\label{sec discrete laplacian}
\vskip .1 in

To have a well-defined discrete Laplacian matrix $\lm$ for a metrized
graph $\ga$, we first choose a vertex set $\vv{\ga}$ for $\ga$ in
such a way that there are no self-loops, and no multiple edges
connecting any two vertices. This can be done for any graph by
enlarging the vertex set by inserting more valence two vertices
whenever needed. We will call such a vertex set $\vv{\ga}$ \textit{optimal}.
If distinct vertices $p$ and $q$ are the end points of an edge, we
call them adjacent vertices.

Note that a weighted graph corresponding to a metrized graph with optimal vertex set
does not have self loops and multiple edges. Such graphs are usually called
simple graphs in the literature.

Let $\ga$ be a graph with $e$ edges and with an optimal vertex set
$\vv{\ga}$ containing $n$ vertices. Fix an ordering of the vertices
in $\vv{\ga}$. Let $\{L_1, L_2, \cdots, L_e\}$ be a labeling of the
edge lengths. The $n \times n$ matrix $\am=(a_{pq})$ given by
\[
a_{pq}=\begin{cases} 0 & \quad \text{if $p = q$, or $p$ and $q$ are
not adjacent}.\\
\frac{1}{L_k} & \quad \text{if $p \not= q$, and $p$ and $q$ are
connected
by} \text{ an edge of length $L_k$}\\
\end{cases}.
\]
is called the adjacency matrix of $\ga$. Let $\dm=\diag(d_{pp})$ be
the $n \times n$ diagonal matrix given by $d_{pp}=\sum_{s \in
\vv{\ga}}a_{ps}$. Then $\lm:=\dm-\am$ is defined to be the \textit{discrete
Laplacian matrix} of $\ga$. That is, $ \lm =(l_{pq})$ where
\[
\displaystyle l_{pq}=\begin{cases} 0 & \; \, \text{if $p \not= q$, and $p$ and $q$
are not adjacent}.\\
-\frac{1}{L_k} & \; \, \text{if $p \not= q$, and $p$ and $q$ are
connected by} \text{ an edge of length $L_k$}\\
\displaystyle -\sum_{s \in \vv{\ga}-\{p\}}l_{ps} & \; \, \text{if $p=q$}
\end{cases}.
\]

The discrete Laplacian matrix is also known as the generalized (or
the weighted) Laplacian matrix in the literature.

Throughout this paper, all matrices will have entries in $\RR$. Given a matrix $\mm$, let $\mm^{T}$, $\mathrm{tr}(\mm)$,
be the transpose and trace
of  $\mm$, respectively.
Let $\idm_n$ be the $n \times n$ identity matrix, and let $\om$ be
the zero matrix (with the appropriate size if it is not specified).

A matrix $\mm$ is called doubly centered, if both row and column
sums are $0$. That is, $\mm$ is doubly centered iff $\mm \ay= \om$
and $\ay^T \mm =\om$, where $\ay=[1,1,\cdots,1]^T$.
\begin{remark}\label{remex disc0}
For any graph $\ga$, the discrete Laplacian matrix $\lm$ is
symmetric and doubly centered. That is, $ \displaystyle \sum_{p \in \vv{\ga}} \lpq
= 0$, for each $q \in \vv{\ga}$, and $\lpq =l_{qp}$ for each $p$,
$q$ $\in \vv{\ga}$.
\end{remark}
In our case, $\ga$ is connected by definition. Thus, $\lm$ is of rank $n-1$.
Although $\lm$ is not invertible, it has generalized inverses. In
particular, it has the \textit{pseudo inverse} $\plm$, also known as the
\textit{Moore-Penrose generalized inverse}, which is uniquely determined by
the following properties:
$$i)  \, \, \lm \plm \lm  = \lm, \quad \,
 ii) \, \, \plm \lm \plm  = \plm, \quad \,
 iii) \, \, (\lm \plm)^{T}   = \lm \plm, \quad \,
 iv) \, \, (\plm \lm)^{T}   = \plm \lm .  $$
The following properties hold for both $\lm$ and $\plm$:
\begin{align*}\label{eqn disc4}
i)  \quad &\lm \text{ and } \plm \text{ are symmetric}, & \qquad \qquad
iii) \quad &\lm \text{ and } \plm \text{ are EP matrices},
\\ ii) \quad &\lm \text{ and } \plm \text{ are doubly centered}, & \qquad \qquad
iv) \quad &\lm \text{ and } \plm \text{ are positive
semi-definite}.
\end{align*}
For more information about  $\lm$ and $\plm$, see the article \cite{C3} and references therein.

I. Gutman and W. Xiao \cite[Lemma 3]{I-W} obtained the
following result when $\lm$ arises from a graph where each edge
length is $1$.
\begin{lemma}\cite[Equation 2.9]{D-M}\label{lem disc1}
Let $\jm$ be an $n \times n$ matrix having each entries $1$, and let $\lm$ be the
discrete Laplacian of a graph (not necessarily with equal edge lengths).
Then $\lm \plm = \plm \lm = \idm -\frac{1}{n}\jm$.
\end{lemma}
As an immediate consequence of \lemref{lem disc1}, we have
\begin{corollary}\label{cor disc1}
Let $\ga$ be a graph and let $\lm$ be the corresponding discrete
Laplacian matrix of size $n \times n$. Then for any $p$, $q$ $\in
\vv{\ga}$,
$ \displaystyle \sum_{s \in \vv{\ga}}
l_{ps}^+ l_{sq}=\begin{cases} -\frac{1}{n} & \text{ if $p \not =q$;}\\
\frac{n-1}{n} & \text{ if $p=q$ .} \end{cases}$
\end{corollary}

The resistance distance has been studied in chemical literature and in
computer science. See the articles \cite{D-M} and \cite{RB1}.
The following fact plays an important role for the
rest of the paper:
\begin{lemma} \cite{RB2} \cite{RB1} \label{lem disc2}
Let $\ga$ be a graph with discrete Laplacian $\lm$ and the
resistance function $r(x,y)$, and let $\hm$ be a generalized inverse
of $\lm$, i.e., $\lm \hm \lm = \lm$. Then
$r(p,q)=\hm_{pp}-\hm_{pq}-\hm_{qp}+\hm_{qq}, \quad \text{for any $p$, $q$ $\in \vv{\ga}$}.$
In particular, for the pseudo inverse $\plm$ we have
$r(p,q)=l_{pp}^+-2l_{pq}^+ + l_{qq}^+, \quad \text{for any $p$, $q$ $\in \vv{\ga}$}.$
\end{lemma}
Note that \lemref{lem disc2} for the pseudo inverse $\plm$ follows from \cite[Theorem A]{D-M}. Similarly, Babi\'{c} {\em et al} \cite{BKLNT}
has \lemref{lem disc2} for $\plm$ and the graphs with all edges of length $1$.

Now, we can easily derive the following analogous result about the
voltage function:
\begin{lemma}\label{lem voltage1}
Let $\ga$ be a graph with the discrete Laplacian $\lm$ and the voltage
function $j_x(y,z)$. Then for any $p$, $q$, $s$ in $\vv{\ga}$,
$j_p(q,s)=l_{pp}^+ - l_{pq}^+ -l_{ps}^+ +l_{qs}^+.$
\end{lemma}
\begin{proof}
By \eqnref{eqn1.1b}, $2j_p(q,s)=r(q,p)+r(s,p)-r(q,s)$. Then the
result follows from \lemref{lem disc2}.
\end{proof}


\section{Identities analogous to the extended Foster's identities}\label{sec Foster Identities}

Let $\ga$ be a metrized graph with an optimal vertex set $\vv{\ga}$ containing $n$ vertices as before.
%
Recall that if we consider $\ga$ as a resistive electric circuit, the resistance of an edge is given by its length.
If $\li$ is the length of an edge $e_i$ with end points $p$ and $q$, then we have
$C_{p q}=C_{q p}=\frac{1}{\li}$ as the conductance of the edge $e_i$.
We write $q \sim p$ when the vertex $q$ is adjacent to the vertex
$p$, i.e., when $p$ and $q$ are the end points of an edge.
We set $C_{pp}=0$ and $C_{pq}=0$ if $p$ and $q$ are not adjacent, and define
$C_p:=\displaystyle \sum_{q \in \vv{\ga}, \, q \sim p} C_{pq}.$
It follows from the definitions that  $C_p=l_{pp}$ and $C_{pq}=-l_{pq}$ if $p \neq q$,
where $\lm=(l_{pq})$ is the discrete Laplacian of $\ga$.

For any $s \in \vv{\ga}$, we \cite[Section 5.5]{C1} proved the following identities by using \lemref{lem voltage1} and the properties of $\lm$ and $\plm$:
\begin{equation*}\label{}
\begin{split}
&\sum_{\substack{p, \, q \in \vv{\ga}\\p \sim q, \, \, \, p<q}}j_p(q,s)C_{pq}=\frac{v-1}{2},
\\&
\sum_{\substack{p, \, q, \, w \in \vv{\ga} \\ w \sim p \sim q, \, w < q}}j_{w}(q,s)
\frac{C_{w p}C_{p q}}{C_{p}}=\frac{v-2}{2},
\\&
\sum_{\substack{p, \, q, \, t, \, w \in \vv{\ga} \\ w \sim p \sim q \sim t, \, w < t}}j_{w}(t,s)
\frac{C_{w p}C_{pq}C_{q t}}{C_{p}C_q}=\frac{v-3}{2}
+\frac{1}{2}\sum_{\substack{p, \, q \in \vv{\ga} \\ p \sim q}} \frac{C_{pq}^2}{C_{p} C_{q}},
\\ &
\sum_{\substack{p, \, q, \, t, \, u, \, w \in \vv{\ga} \\ u \sim p \sim w \sim q \sim t \\ u < t}}j_{u}(t,s)
\frac{C_{u p}C_{p w} C_{wq}C_{q t}}{C_{p} C_w C_q}=\frac{v-4}{2}
+\frac{1}{2}\sum_{\substack{p, \, q \in \vv{\ga} \\ p \sim q}} \frac{C_{pq}^2}{C_{p} C_{q}}
+\frac{1}{2}\sum_{\substack{p, \, q, \, w \in \vv{\ga} \\ p \sim q \sim w \sim p}} \frac{C_{pw}C_{wq}C_{qp}}{C_{p} C_w C_{q}}.
\end{split}
\end{equation*}
Note that the formulas above are for graph paths containing at most $4$ edges. A certain pattern involving conductances can be seen
in these formulas.
In the rest of the paper, we will extend formulas above, and obtain a new proof of the extended Foster's identities as a corollary.

Let $\mathrm{P}=(\mathrm{p}_{it})$ be the $n \times n$ matrix given by $\mathrm{p}_{it}=\displaystyle \frac{C_{it}}{C_i}$, and
let $\mathrm{P^k}=(\mathrm{p}_{it}^{(k)})$ be the $k$-th power of $\mathrm{P}$. It follows from the definitions that
$\mathrm{p}_{it}^{(k+1)}= \displaystyle \sum_{m_{1}, \ldots, m_{k}=1}^n\frac{C_{i m_1} C_{m_1 m_2} \cdots \, C_{m_{k} t}}{C_{i} C_{m_1} \cdots \, C_{m_k}}$ whenever $k \geq 1$. Moreover, we have $\displaystyle \sum_{t=1}^n\mathrm{p}_{it}^{(k)}=1$  and that $ \displaystyle \sum_{i=1}^n C_i \mathrm{p}_{it}^{(k)}=C_{t}$ for any $k \geq 1$. As it is also given in the article \cite{BCEG}, $\mathrm{P}$ is known as \textit{transition probability matrix} of the reversible Markov chain associated to $\ga$ when $\ga$ is considered as an electrical network, and $\mathrm{P^k}$ is known as the \textit{k-step transition probability matrix}. The value $\mathrm{p}_{it}^{(k)}$ is the probability that the Markov chain attains vertex $t$ at $k$-th step after starting from vertex $i$.
\begin{proposition}\label{prop trans matrix}
For any $k \geq 1$, we have
$\displaystyle \sum_{i, \, t=1}^n C_i l_{it}^+ \mathrm{p}_{it}^{(k+1)}=1- \mathrm{tr} (\mathrm{P^k}) + \sum_{i, \, t=1}^n C_i l_{it}^+ \mathrm{p}_{it}^{(k)}.$
\end{proposition}
\begin{proof}
We first note that
\begin{equation}\label{eqn trans1}
\begin{split}
\sum_{i=1}^n l_{it}^+ C_{iq}& =l_{q t}^+l_{q q}-\sum_{i=1}^n l_{it}^+ l_{iq}
=C_{q} l_{q t}^+ -\begin{cases} -\frac{1}{n} & \text{ if $t \not =q$;}\\
\frac{n-1}{n} & \text{ if $t=q$.} \end{cases}
\end{split}
\end{equation}
The last equality follows from \corref{cor disc1}. For any $k \geq 1$, we have
\begin{equation*}\label{}
\begin{split}
\sum_{i,\, t=1}^n C_i l_{it}^+ \mathrm{p}_{it}^{(k+1)} & =\sum_{t, \, m_{1}, \ldots, m_{k}=1}^n\frac{C_{m_1 m_2} \cdots \, C_{m_{k} t}}{ C_{m_1} \cdots \, C_{m_k}} \sum_{i=1}^n l_{it}^+ C_{i m_1}, \quad \, \text{Then by using \eqnref{eqn trans1},}
\\ &=\sum_{i, \, t=1}^n C_i l_{it}^+ \mathrm{p}_{it}^{(k)}+\frac{1}{n} \sum_{m_{1}, \ldots, m_{k}=1}^n \sum_{\substack{t=1 \\ t\neq m_1}}^n \frac{C_{m_1 m_2} \cdots \, C_{m_{k} t}}{ C_{m_1} \cdots \, C_{m_k}}-\frac{n-1}{n} \sum_{m_{1}=1}^n \mathrm{p}_{m_1 m_1}^{(k)}
\\ &=\sum_{i, \, t=1}^n C_i l_{it}^+ \mathrm{p}_{it}^{(k)}+\frac{1}{n} \sum_{t, \, m_{1}, \ldots, m_{k}=1}^n\frac{C_{m_1 m_2} \cdots \, C_{m_{k} t}}{ C_{m_1} \cdots \, C_{m_k}}
-\tr (\mathrm{P^k})
\\ &=\sum_{i, \, t=1}^n C_i l_{it}^+ \mathrm{p}_{it}^{(k)}+\frac{1}{n} \sum_{t, \,  m_{1}=1}^n \mathrm{p}_{m_1 t}^{(k)} -\tr (\mathrm{P^k}).
\end{split}
\end{equation*}
Then the result follows from the fact that $\displaystyle \sum_{t, \, m_{1}=1}^n \mathrm{p}_{m_1 t}^{(k)}=\sum_{m_1=1}^n 1=n$ for any $k \geq 1$.
\end{proof}
We have $\displaystyle \sum_{i,\, t=1}^n C_i l_{it}^+ \mathrm{p}_{it}=\sum_{i, \, t=1}^n l_{it}^+ C_{it}=-\sum_{i, \, t=1}^n l_{it}^+ l_{it}+\sum_{i}^n C_{i} l_{ii}^+=-n+1+\sum_{i}^n C_{i} l_{ii}^+$, where the last equality follows from \corref{cor disc1}.
Thus, the successive application of \propref{prop trans matrix} gives the following equality for any $k \geq 1$:
\begin{equation}\label{eqn trans2}
\begin{split}
\sum_{i, \, t=1}^n C_i l_{it}^+ \mathrm{p}_{it}^{(k)}=k-n - \sum_{i=1}^{k-1}\tr (\mathrm{P^{i}})+\sum_{i=1}^n C_{i} l_{ii}^+.
\end{split}
\end{equation}
Our main result is the following theorem which gives an identity analogues to the extended Foster's identities. Unlike extended Foster's identities which involve the resistance values, our formula involves the voltage values.
\begin{theorem}\label{thm main}
Let $s \in \vv{\ga}$. For any $k \geq 1$, we have
$ \displaystyle \sum_{i, \, t=1}^n C_{i} j_i(s,t) \mathrm{p}_{it}^{(k)}=n-k+\sum_{i=1}^{k-1}$ $\tr (\mathrm{P^i}).$
\end{theorem}
\begin{proof}
Let $s$ be a given vertex in $\vv{\ga}$. For each $k \geq 1$, we have
\begin{equation*}\label{}
\begin{split}
\sum_{i, \, t=1}^n C_{i} j_i(s,t) \mathrm{p}_{it}^{(k)} &= \sum_{i, \, t=1}^n C_{i} \big(l_{ii}^+ -l_{it}^+ -l_{is}^+ +l_{ts}^+ \big) \mathrm{p}_{it}^{(k)}
, \quad \text{by \lemref{lem voltage1}}.
\\ &=\sum_{i=1}^n C_{i} (l_{ii}^+ -l_{is}^+) -\sum_{i, \, t=1}^n C_{i} (l_{it}^+ -l_{ts}^+) \mathrm{p}_{it}^{(k)}
, \quad \text{since $\displaystyle \sum_{t=1}^n \mathrm{p}_{it}^{(k)}=1$.}
\\ &= \sum_{i=1}^n C_{i} l_{ii}^+ -\sum_{i, \, t=1}^n C_{i} l_{it}^+ \mathrm{p}_{it}^{(k)},
\quad \text{since $\displaystyle \sum_{i=1}^n C_{i} \mathrm{p}_{it}^{(k)}=C_{t}$.}
\end{split}
\end{equation*}
Then the result follows from \eqnref{eqn trans2}.
\end{proof}
Next, we will state the extended Foster's identities as a corollary of \thmref{thm main}.
\begin{corollary}\label{cor main}
For any $k \geq 1$, we have
$\displaystyle \frac{1}{2}\sum_{i, \, t=1}^n C_{i} r(i,t) \mathrm{p}_{it}^{(k)}=n-k+\sum_{i=1}^{k-1} \tr (\mathrm{P^i}).$
\end{corollary}
\begin{proof}
By \eqnref{eqn1.1}, we have $r(i,t)= j_i(t,s) + j_t(i,s)$ for any $s \in \ga$. Thus, the result follows from \thmref{thm main}.
\end{proof}
Note that the formula given in \corref{cor main} was originally proved by E. Bendito, A. Carmona, A. M. Encinas and J. M. Gesto \cite[Proposition 2.3]{BCEG}. Next, we will use their methods to provide a second proof for \thmref{thm main}. First, we recall some related definitions and results.

Let \textbf{e} be the $n \times 1$ vector whose entries equal to $1$, and let $\textbf{e}^i$ denote the $n \times 1$
$i$-th unit vector for each $i=1, \ldots , n$. That is, $\textbf{e}^i_k=1$ if $k=i$ and $\textbf{e}^i_k=0$ if $k \neq i$.
Suppose $\nu^i$ is the unique solution of the linear system $\lm u= \textbf{e}-n \textbf{e}^i$ satisfying $\nu^i_i=0$.
The solution $\nu^i$ is called \cite{BCEG} the \textit{equilibrium measure} of the set $\vv{\ga} \backslash \{ i \}$.

Since $r(i,t)=\frac{1}{n} (\nu^i_t+\nu^t_i)$ \cite[Proposition 2.1]{BCEG}, we obtain the following proposition by using
\eqnref{eqn1.1b}:
\begin{proposition}\label{propeqn voltage}
For any $i$, $j$, $t$ in $\vv{\ga}$, we have
$$2 j_i(s,t)=\frac{1}{n} \big(\nu^i_s + \nu^s_i + \nu^i_t + \nu^t_i - \nu^s_t - \nu^t_s \big).$$
\end{proposition}
On the other hand, we have \cite[Proof of Proposition 2.3]{BCEG} that
\begin{equation}\label{eqn bceg}
\begin{split}
\frac{1}{n} \sum_{i, \, t=1}^n C_i \nu^i_t \mathrm{p}_{it}^{(k)} =n-k + \sum_{i=1}^{k-1} \tr (\mathrm{P^i}).
\end{split}
\end{equation}
The second proof of \thmref{thm main} will be as follows.
For any $k \geq 1$, we have
\begin{equation*}
\begin{split}
\sum_{i, \, t=1}^n C_{i} j_i(s,t) \mathrm{p}_{it}^{(k)} &= \frac{1}{2 n} \sum_{i, \, t=1}^n C_{i} \big(\nu^i_s + \nu^s_i + \nu^i_t + \nu^t_i - \nu^s_t - \nu^t_s \big) \mathrm{p}_{it}^{(k)}
, \quad \text{by \propref{propeqn voltage}}.
\\&=\frac{1}{n} \sum_{i, \, t=1}^n C_{i} \nu^i_t \mathrm{p}_{it}^{(k)}
, \quad \text{since $C_i \mathrm{p}_{it}^{(k)}=C_t \mathrm{p}_{ti}^{(k)}$ and $\sum_{t=1}^n \mathrm{p}_{it}^{(k)} =1$}.
\\ &=n-k + \sum_{i=1}^{k-1} \tr (\mathrm{P^i}), \quad \text{by \eqnref{eqn bceg}}.
\end{split}
\end{equation*}



\end{document}